\documentclass[11pt]{amsart}
\usepackage{epsfig,amsmath,amssymb,amsthm,epsf}
\usepackage{graphicx,array}

\newcommand{\zz}[1]{\mathbb #1}
\newtheorem{proposition}{Proposition}
\newtheorem{lemma}{Lemma}
\newtheorem{conjecture}{Conjecture}
\newtheorem{theorem}{Theorem}
\newtheorem{claim}{Claim}

\begin{document}
\title{{A two-species competition model on $\mathbb{Z}^d$ }}
\author{{George Kordzakhia}}
\address{University of California\\
Department of Statistics\\
Berkeley CA }
\email{kordzakh@stat.berkeley.edu}
\author{Steven P. Lalley} 
\address{University of Chicago\\ Department of Statistics \\ 5734
University Avenue \\
Chicago IL 60637}
\email{lalley@galton.uchicago.edu}
\date{\today}
\maketitle
\begin{abstract}
We consider a two-type stochastic competition model on the integer lattice  $\mathbb{Z}^d$.
 The model describes the space evolution of two ``species'' competing for territory along their boundaries. Each site of the space may contain only one representative (also referred to as a particle)  of either type. The spread mechanism for both species is the same: each particle produces offspring  independently of other particles  and can place them only at the neighboring sites that are either unoccupied,  or occupied by particles of the opposite  type. In the second case, the old particle is killed by the newborn. The rate of birth for each particle is equal to the number of neighboring sites  available for expansion.
The main  problem we address  concerns the possibility of the  long-term coexistence of the two species. 
We have shown  that if we start the process  with  finitely  many   representatives  of each type, then, 
under the assumption that the limit set in the corresponding first passage percolation model is  uniformly curved,
there is  positive probability of coexistence. 
\end{abstract}


\section{Introduction}\label{sec:intro}
The \emph{voter model} is an interacting
particle system in which individuals (particles) of two species,
\emph{Red} and \emph{Blue}, compete for ``territory'' on a (locally
finite) graph. At each time $t\geq 0$, every vertex (site) of the
graph is occupied by a single particle, either Red or Blue. At any
time, a particle of color $i$ at a vertex $x$ may spontaneously die,
at rate equal to the degree of $x$, and be replaced by a clone of a
randomly chosen neighbor. Thus, a vertex of color $i$ spontaneously
flips to the opposite color $j$ at rate equal to the number of
neighboring vertices of color $j$. See \cite{liggett} for a formal
construction of this process and an exposition of its basic
properties.

The \emph{Richardson model} \cite{richard} was introduced as a model
for the spatial spread of a population in a favorable environment. The
environment is once again a locally finite graph. At any time a vertex
may be occupied by \emph{at most} one particle (some vertices may be
unoccupied); all particles are of the same species.  Once occupied, a
vertex remains occupied forever. Each unoccupied vertex is
spontaneously occupied at instantaneous rate equal to the number of
occupied neighbors.

In this paper we study a hybrid of the voter and Richardson models on
the integer lattice $\zz{Z}^{d}$, which we dub the \emph{two-species
competition model}, or simply the \emph{competition model}. The
dynamics are as in the voter model, but unlike the voter model,
vertices may be unoccupied. An unoccupied vertex is colonized at rate
equal to the number of occupied neighbors, as in the Richardson model;
at the instant of first colonization, the vertex flips to the color of
a randomly chosen occupied neighbor. Once occupied, a vertex remains
occupied forever, but its color may flip, as in the voter model: the
flip rate is equal to the number of neighbors occupied by particles of
the opposite color.  The state of the system at any time $t$ is
given by the pair $(R (t),B (t))$, where $R (t)$ and $B (t)$ denote
the set of sites occupied by Red and Blue particles,
respectively. Note that the set $R (t)\cup B (t)$ of occupied sites
evolves precisely as in the Richardson model, and so the growth of
this set is governed by the same \emph{Shape Theorem} (see section
\ref{sec:richardson} below) as is the Richardson model.

Our primary interest is in the possibility of long-term coexistence of
the two species, given initial conditions in which only finitely many
vertices are occupied (with at least one vertex of each color). It is
clear that at least one of the two species must survive, and that for
any nondegenerate finite initial configuration of colors there is
positive probability that Red survives and positive probability that
Blue survives. However, it is not at all obvious (except perhaps in
the case where the ambient graph on which the competition takes place
is the integer lattice $\zz{Z}$ -- see section \ref{sec:1D} below)
that the event of mutual survival has positive probability.  Our main
result concerns the competition model on the graph $\zz{Z}^{d}$. Say
that a compact, convex set $\mathcal{S}$ with boundary $\partial
\mathcal{S}$ is \emph{uniformly curved} if there exists $\varrho
<\infty$ such that for every point $z\in \partial \mathcal{S}$ there
is a ball of radius $\varrho$ with $z$ on its surface that contains
$\mathcal{S}$.

\begin{theorem}  \label{Newman:Theorem1}
If the limit shape $\mathcal{S}$ for the Richardson model is uniformly
curved, then for any nondegenerate initial finite configuration the
event of mutual survival of the two species has positive probability.
\end{theorem}

The proof will be carried out in sections
\ref{sec:preliminaries}--\ref{sec:proof} below.  Theorem
\ref{Newman:Theorem1} is by no means a complete solution to the
coexistence problem, because it remains unknown whether the limit
shape $\mathcal{S}$ for the Richardson model is uniformly curved, or
even if its boundary $\partial \mathcal{S}$ is strictly
convex. Nevertheless, simulations give every indication that it is,
and \cite{lalley} suggests a possible explanation of what lies behind
the strict convexity of $\partial \mathcal{S}$.

The two-species complete model is superficially similar to the
\emph{two-type Richardson model} studied by Haggstrom and Pemantle
\cite{haggstrom-pemantle}, but differs in that it allows displacement
of colors on occupied sites: in the two-type Richardson model, once a
vertex is occupied by a particle (either red or blue) it remains
occupied by that color forever. The main result of
\cite{haggstrom-pemantle} is similar to Theorem \ref{Newman:Theorem1},
but requires no hypothesis about the Richardson shape: it states that
mutual unbounded growth has positive probability. Because no
displacements are allowed, the behavior of the two-type Richardson
model is very closely tied up with the first-passage percolation
process with exponential edge passage times. The two-species
competition model is also closely related to first-passage
percolation, but the connection is less direct, because the
possibility of displacements implies that not only the first passages
across edges play a role in the evolution.

\begin{figure}
\caption{Competition Model with $R(0)=\{ (0,0) \},\ B(0)=\{ (1,0) \}$.
\label{800competition}} \vspace{1.5ex}
\begin{tabular}{cc}
\includegraphics[width=2.5in, height=2.5in]{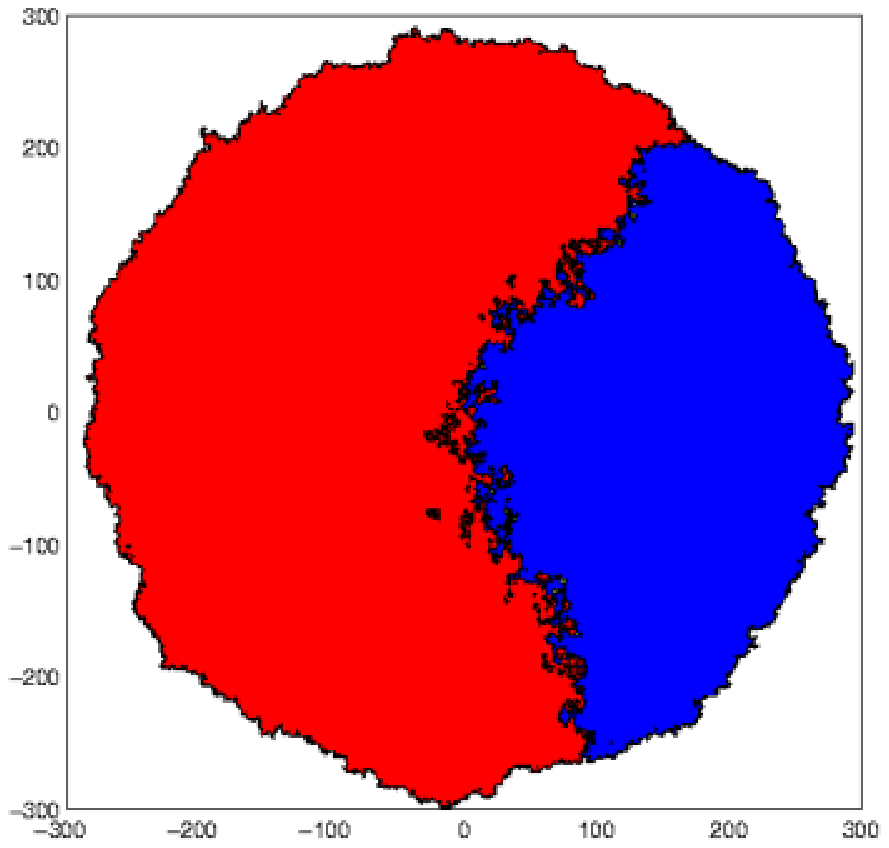}
\includegraphics[width=2.5in, height=2.5in]{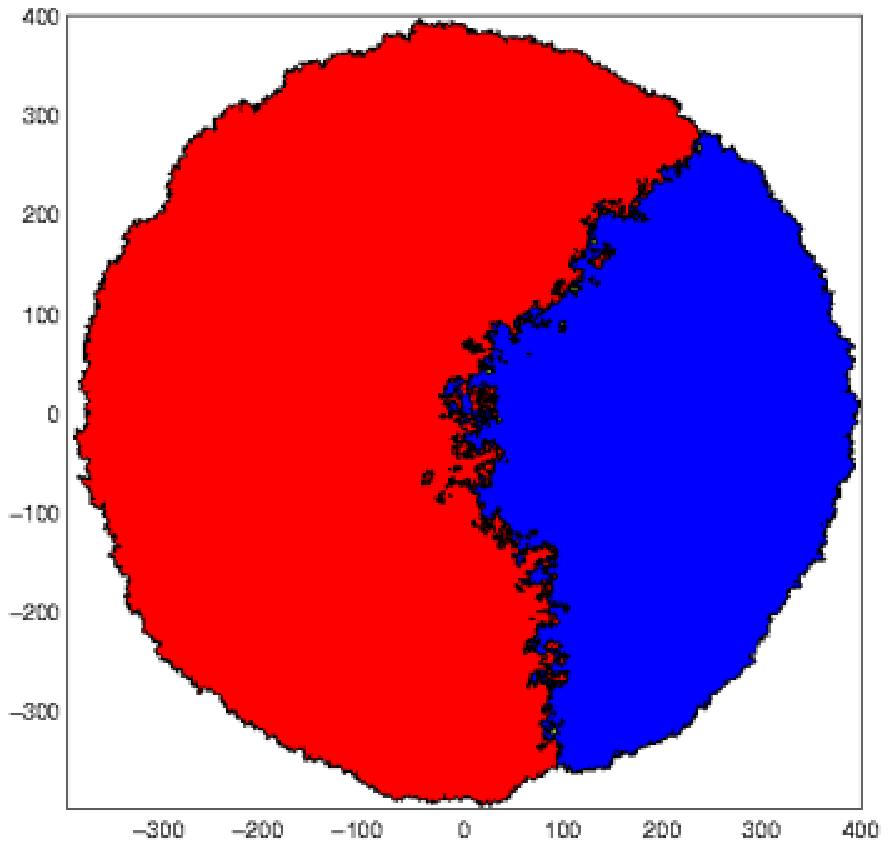} 
\end{tabular}
\end{figure}

We have run simulations of the competition model with initial
configuration $R (0) =\{(0,0) \}$ and $B (0)=\{(1,0) \}$. Figure
\ref{800competition} shows two snap shots of the same realization of
the process taken at the times when the region occupied by both types,
$R(t) \cup B(t)$, hits the boundary of the rectangles $[-300,\ 300]
\times [-300,\ 300]$ and $[-400,\ 400] \times [-400,\ 400]$
respectively. Observe that the overall shape of the red and the blue
clusters did not change significantly. We believe that the shape of
the regions occupied by the red and blue types stabilizes as times
goes to infinity.

\begin{figure}
\caption{Competition Model with $R(0)=\{ (-2,-2),\ (2,2) \}$ and
$B(0)=\{ (-2,2),\ (2,-2) \}$.  \label{Mixedcompetition}}
\vspace{1.5ex}
\begin{tabular}{cc}
\includegraphics[width=2.4in, height=2.4in]{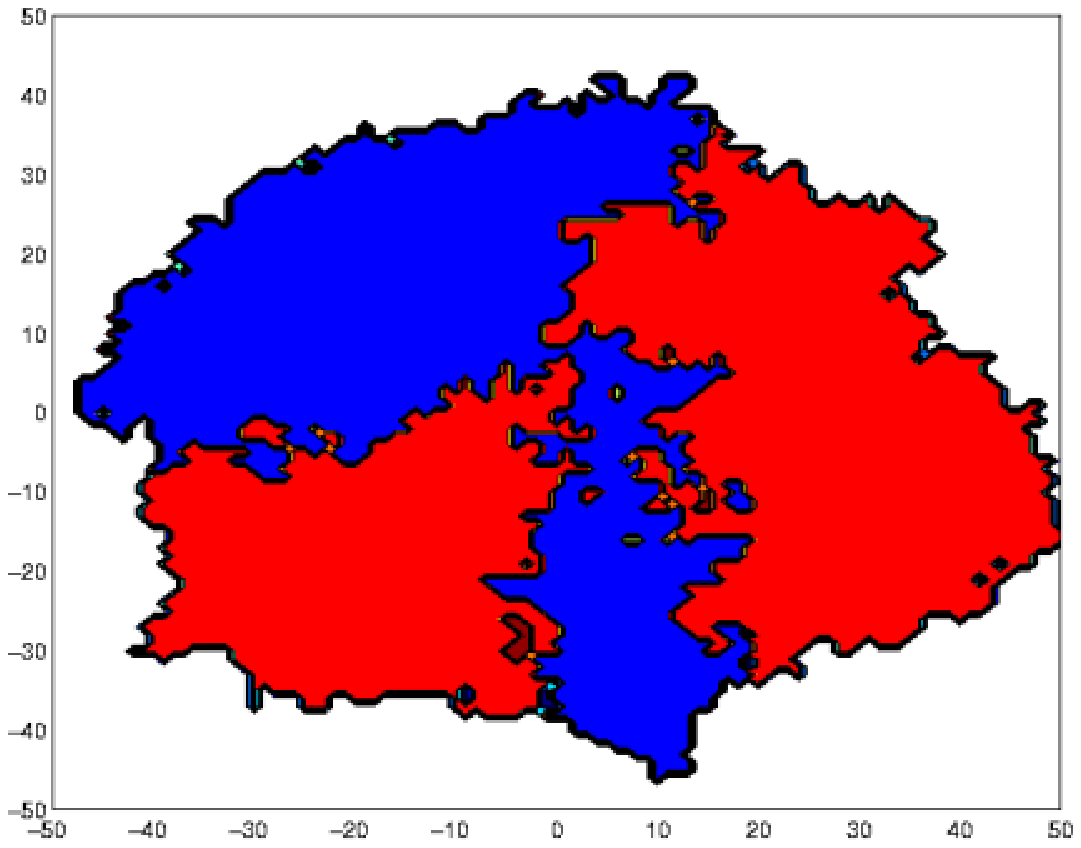} 
 \includegraphics[width=2.4in, height=2.4in]{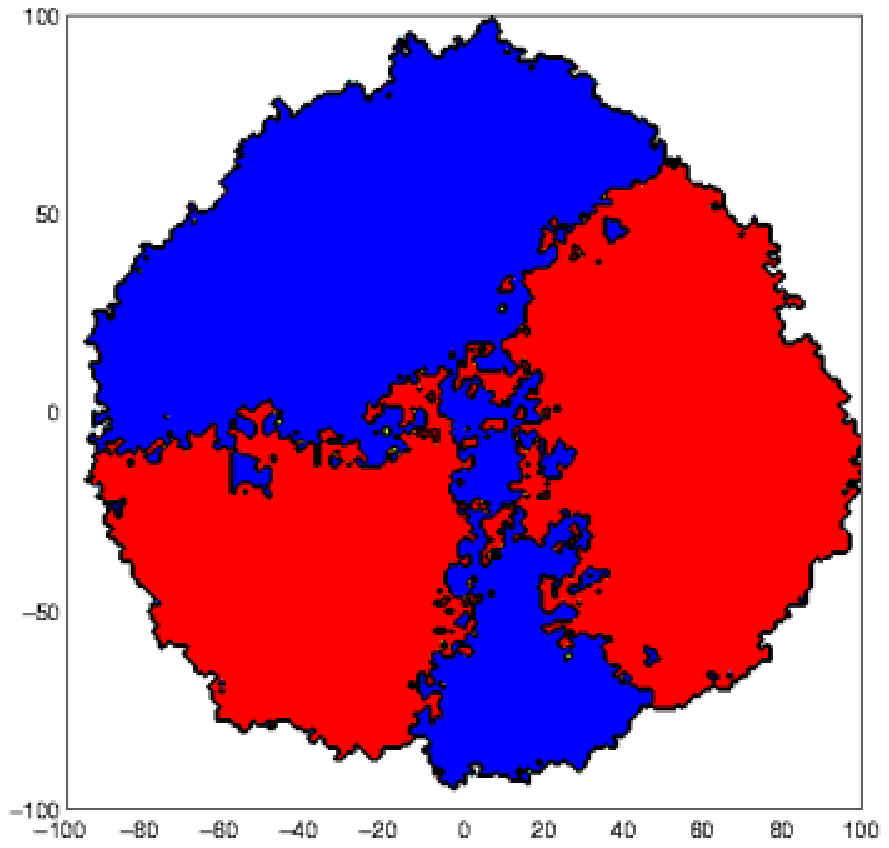}  \\
 \includegraphics[width=2.4in, height=2.4in]{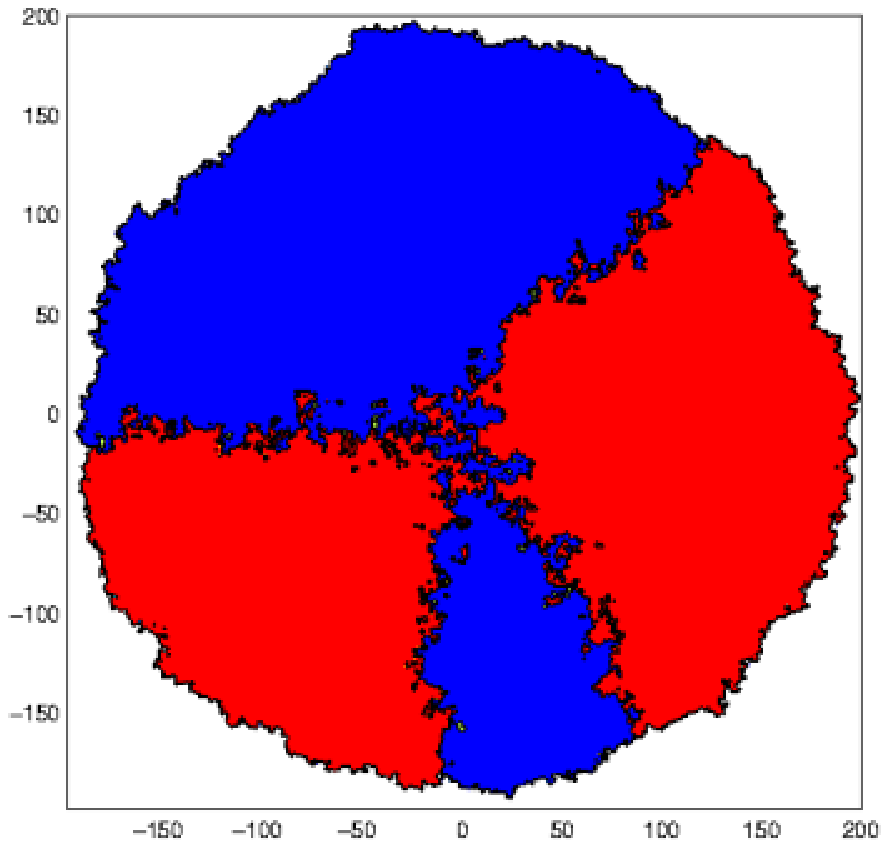}
\end{tabular}x
%\parbox{6in}{Note: }
\end{figure}

For any subset $Z\subset \zz{R}^{d}$, define
\[
        \hat{Z}=\{x\in \zz{R}^{d}\, : \, \text{dist} (x,Z)\leq 1/2\},
\]
where $\text{dist}$ denotes distance in the $L^{\infty}-$norm on
$\zz{R}^{d}$. For any subset $Z\subset \zz{R}^{d}$ and any scalar
$s>0$, let $Z/s=\{y/s \, : \, y\in Z\}$. 
\begin{conjecture}
There exist random sets $\tilde{R}$ and $\tilde{B}$ such that with
probability one
\begin{align}\label{conj}
        \lim_{t \rightarrow \infty}\hat{R}(t)/t&= \tilde{R}, \\
        \lim_{t \rightarrow \infty} \hat{B}(t)/t&= \tilde{B}, \quad \text{and}\\
        \tilde{R}\cup  \tilde{B}&= \mathcal{S}. 
\end{align}
\end{conjecture}
If this is true, we expect that the limit sets $\tilde{R}$ and
$\tilde{B}$ will be finite unions of angular sectors, as the
simulation results shown in Figure \ref{800competition} suggest. The
sizes and directions of these angular sectors (and even their number)
will, we expect, be random, with distributions depending on the
initial configuration. This is illustrated by simulation results
summarized in Figure \ref{Mixedcompetition} with initial configuration
\begin{align*}
        R (0)&=\{  (-2,-2),\ (2,2)\}, \\
        B (0)&=\{ (-2,2),\ (2,-2)\}.
\end{align*}
Three time progressive snap shots of the process were taken. The plots
in the figure \ref{Mixedcompetition} suggest that stabilization of
the shape was taking place on the considered time interval.

\section{The Competition Model on $\zz{Z}^{1}$}\label{sec:1D}

The coexistence problem for the Competition Model in one dimension is
considerably simpler than in higher dimensions. Since the limit shape
of the Richardson model in one dimension is  an interval, no auxiliary
hypothesis is needed.

\begin{proposition}\label{proposition:1D}
For any nondegenerate finite initial configuration on $\zz{Z}$, the
event of mutual survival in the two-species Competition Model has
positive probability.
\end{proposition}

\begin{proof}
Without loss of generality, we may assume that the initial
configuration consists of a finite interval of red sites with
rightmost point $-1$ and a finite interval of blue sites with leftmost
point $0$, since [a translate of] such a configuration may be reached
in finite time, with positive probability, from \emph{any}
nondegenerate initial configuration. Let $X_{t}$ and $Y_{t}$ be the
left- and right-most occupied sites (of either color) at time $t$, and
let $Z_{t}$ be the leftmost blue site. Note that as long as
$X_{t}<Z_{t}<Y_{t}$, there will be both red and blue sites: all sites
to the left of $Z_{t}$ are red, and all sites to the right are
blue. Each of the processes $X_{t}$ and $Y_{t}$ is a pure jump
process, with jumps of size $1$ occurring at rate $1$; hence, with
probability one, as $t \rightarrow \infty$,
\begin{align*}
        X_{t}/t &\longrightarrow -1 \quad \text{and}\\
        Y_{t}/t &\longrightarrow 1.
\end{align*}
The process $Z_{t}$ behaves, up to the time of first exit from
$(X_{t},Y_{t})$, as a continuous-time simple nearest-neighbor random
walk on the integers. Consequently, there is positive probability that
$Z_{t}$ never exits the interval $(X_{t},Y_{t})$. But on this event,
both species survive.
\end{proof}

This simple argument clearly shows what the difficulty in higher
dimensions will be: In one dimension, the interface between
(connected) red and blue clusters is just a point; but in higher
dimensions, it will in general be a hypersurface, whose time evolution
will necessarily be somewhat complicated.

\section{Preliminaries}\label{sec:preliminaries}

\subsection{Graphical Constructions}\label{sec:graphicalConstructions}

The Richardson model, the voter model, and the two-species competition
model all admit \emph{graphical contructions} using \emph{percolation
structures}.  Such constructions make certain comparison arguments and
duality relations transparent. We briefly review the construction
here, primarily to emphasize that the same percolation structure can
be used to simultaneously build versions of all three processes with
all possible initial configurations. See, for instance,
\cite{DurrettBook} for further details in the case of the Richardson
model and the voter model.

The \emph{percolation structure} $\Pi$ is an assignment of
independent, rate-$1$ Poisson processes to the directed edges $xy$ of
the lattice $\zz{Z}^{d}$. (For each pair $\{x,y \}$ of neighboring
vertices, there are two directed edges $xy$ and $yx$.)  Above each
vertex $x$ is drawn a timeline, on which are placed marks at the
occurrence times $T^{xy}_{i}$ of the Poisson processes attached to
directed edges emanating from $x$; at each such mark, an arrow is
drawn from $x$ to $y$. A \emph{directed path} through the percolation
structure $\Pi$ may travel upward, at speed $1$, along any timeline,
and may (but does not have to) jump across any outward-pointing arrow
that it encounters. A \emph{reverse path} is a directed path run
backward in time: thus, it moves downward along timelines and jumps
across inward-pointing arrows. A \emph{voter-admissible} path is a
directed path that does not pass any inward-pointing arrows. Observe
that for each vertex $y$ and each time $t>0$ there is a unique
voter-admissible  path beginning at time $0$ and terminating
at $(y,t)$: its reverse path is gotten by traveling
downward along timelines, starting at $(y,t)$, jumping across all
inward-pointing arrows encountered along the way. 

\noindent 
\emph{Richardson Model:} A version $Z (t)$ of the Richardson model
with initial configuration $Z (0)=\zeta $ is obtained by setting
$Z(t)$ to be the set of all vertices $y$ such that there is a directed
path in the percolation structure $\Pi$ that begins at $(x,0)$ for
some $x\in \zeta$ and ends at $(y,t)$.

\noindent 
\emph{Voter Model:} A version $(R (t), B (t))$ of the voter model
with initial configuration $R (0)=\beta $, $B (0)=\beta^{c}$ is gotten
by defining $B (t)=R (t)^{c}$ and  $R (t)$ to be the set of all
vertices $y$ such that the unique voter-admissible path terminating at
$(y,t)$ begins at  $(x,0)$ for some  $x\in \beta$.

\noindent \emph{Two-Species Competition Model:} Fix an initial
configuration $R(0)=\xi $, $B (0)=\zeta$. Erase all arrows that lie
\emph{only} on  paths that begin at points $(x,0)$ such
that $x\not \in \xi \cup \zeta$; denote the resulting sub-percolation
structure $\Pi_{\xi ,\zeta}$. Define $R (t)$ (respectively, $B (t)$)
to be the set of all vertices $y$ such that there is a
voter-admissible path relative to $\Pi_{\xi ,\zeta}$ that ends at
$(y,t)$ and starts at $(x,0)$ with $x\in \xi$ (respectively, $x\in
\zeta $).

The graphical construction yields as by-products comparison principles
for the Richardson, voter, and competition models. First, the set $R
(t)\cup B (t)$ of vertices occupied by either Red or Blue particles at
time $t$ in the competition model coincides with the set $Z (t)$ of
occupied vertices in the Richardson model when $Z(0)=R(0)\cup B
(0)$. Second, if $R^{(1)}(t), B^{(1)}(t)$ is the voter model with
initial configuration $R^{(1)}(0)=\xi^{(1)}$, and $R
(t),B (t)$ is the competition model with initial configuration $R
(0)=\xi \supset \xi^{(1)}$, $B (0)=\zeta \subset \xi^c$, then for all $t>0$,
\[
        R (t)\supset R^{(1)}(t).
\]

\subsection{Voter Model: Invasion Times}\label{sec:voter}

How long does it take for a Red vertex to be overrun by Blue?
Clearly, for either the voter model or the competition model the
answer will depend, at least in part, on how far away the nearest Blue
vertices are. The comparison principle implies that, for any given
value $\rho$ of the distance to the nearest blue vertex, the worst
case (for either model) is the voter model with initial configuration
$R (0)=D (0,\rho)$ and $B(0)=D (0,\rho)^{c}$, where $D (x,\rho)$
denotes the disk of radius $\rho$ centered at $x$ (more precisely, its
intersection with the lattice $\zz{Z}^{d}$).

\begin{lemma}\label{lemma:invasion}
Fix $\beta \in (\frac{1}{2} ,1)$, and denote by $(R (t),B (t))$ the
state of the voter model at time $t$. There exist constants
$c_{1},c_{2}>0$ (depending on $\beta$) such that for all $\rho \geq 1$
and all $t\in [0,\rho]$, if $R (0)$ contains the disk of radius
$\rho^{\beta}$ centered at the vertex $x$, then
\[
        P\{x\not \in R (t) \}\leq c_{1}
                \exp\{-c_{2}\rho ^{\beta - 1/2 } \}.
\]
\end{lemma}

\noindent \textbf{Remark.} This holds for \emph{any} norm on
$\zz{R}^{d}$, not just the Euclidean norm: in particular, it holds for
the \emph{Richardson norm} defined below. The constant $c_{2}>0$ may,
of course, depend on the norm.

\begin{proof}
The dual process of the voter model is the coalescing random walk (see
\cite{DurrettBook} or \cite{liggett}). Thus, the probability that the
vertex $x$ is blue at time $t$ coincides with the probability that a
continuous-time simple random walker started at $x$ at time $0$ will
land in the set $B (0)$ at time $t$. (This is not difficult to deduce
directly from the graphical construction above: the event $x \in B(t)$
occurs if and only if the reverse voter-admissible path started 
at $(x,t)$ will terminate at $(y,0)$ for some $y\in B (0)$; but the
reverse voter-admissible path is a simple random walk.) Hence, if
$R(0)\supset D (x,\rho^{\beta})$ then this probability is dominated by the
probability that the continuous-time simple random walk  exits the ball
$D(x,\rho^{\beta})$ by time $\rho$. 
\end{proof}

\subsection{Richardson Model: Shape Theorem}\label{sec:richardson} 

The first-order asymptotic behavior of the Richardson model on the
integer lattice $\zz{Z}^{d}$ is described by the \emph{Shape Theorem}
\cite{cox}. Denote by $Z (t)$ the set of vertices of $\zz{Z}^{d}$ that
are occupied at time $t$, and by $P_{\zeta}$ the probability measure
describing the law of the process given the initial condition $Z
(0)=\zeta$. 
For any subset $Z\subset \zz{R}^{d}$, define
\[
        \hat{Z}=\{x\in \zz{R}^{d}\, : \, \text{dist} (x,Z)\leq 1/2\},
\]
where $\text{dist}$ denotes distance in the $L^{\infty}-$norm on
$\zz{R}^{d}$. For any subset $Z\subset \zz{R}^{d}$ and any scalar
$s>0$, let $Z/s=\{y/s \, : \, y\in Z\}$. 

\begin{theorem} [The Shape Theorem] \label{RichShape}
There exists
a~nonrandom compact convex set $\mathcal{S}\subset \mathbb{R}^d$,
invariant under permutation of and reflection in the coordinate
hyperplanes, and with non-empty interior, such that for any finite
initial configuration $Z (0)=\zeta$ and any
$\epsilon>0$, with $P_{\zeta}-$probability one, eventually (i.e., for
all sufficiently large $s$)
\[ 
 (1-\epsilon)\cdot \mathcal{S} 
        \subseteq \frac{\hat{Z}(s)}{s}
        \subseteq (1+\epsilon)\cdot \mathcal{S} .
\]   
\end{theorem}

The exact shape of the limiting
set $\mathcal{S}$ remains unknown. A simple argument shows that
$\mathcal{S}$ is convex, but nobody has succeeded in proving that it is
\emph{strictly} convex.  Let $|\cdot|$ be the norm on $\mathbb{R}^d$ associated
with the shape set $\mathcal{S}$, that is, for $x \in \mathbb{R}^d$,
$|x|=\inf \{ t: x \in \mathcal{S} \cdot t \}$. That this is in fact a
norm follows  from the convexity of $\mathcal{S}$. The Shape Theorem is
equivalent to the statement that the set of occupied sites grows at
speed one, relative to the norm $| \cdot |$, in every direction.

The Richardson model admits a description as a first passage
percolation model, as follows. To each edge of the
lattice $\mathbb{Z}^d$, attach a mean one exponential random variable,
the `` passage time'', in such a way that the passage times of
distinct edges are mutually independent. For any self-avoiding path
$\gamma$, define the traversal time $\tau (\gamma)$ to be the sum of the
passage times of the edges in $\gamma$. For any finite set $\zeta$ of
vertices and any vertex $x$, define the passage time $T(\zeta ,x)$
from $\zeta$ to $x$ to be the infimum of the traversal times 
$\tau (\gamma )$ of all self-avoiding paths connecting $x$ to
$\zeta$. A version of the Richardson model $Z (t)$ with initial
configuration $Z (0)=\zeta$ is given by
\[
        Z (t):= \{x\in \zz{Z}^{d}\, : \, T (\zeta ,x)\leq t\}.
\]
The  first-passage percolation representation gives
simultaneous realizations of Richardson evolutions for all initial
configurations. Since traversal times of paths are the same backwards
and forwards, the following \emph{duality property} is immediate: For
any finite subsets $F,G \subset \zz{Z^{d}}$ and any $t>0$,
\begin{equation}\label{eq:richardsonDuality}
        P_{F}\{ Z(t) \cap G\  =\emptyset \}=
        P_{G}\{ Z(t) \cap F\  =\emptyset \}.
\end{equation}

Kesten \cite{kesten} and Alexander \cite{alexander} have established
large deviation results for the passage times in first passage
percolation that specialize to the Richardson model as follows.  
\begin{theorem} \label{Kesten}
There exist constants $c_1$ and $c_2 >0$ such that for any
$\epsilon>0$, $y \in \mathbb{Z}^d$ and ${|y|}^{1/2+\epsilon}< t<
|y|^{3/2-\epsilon}$,
\[ 
        P \left\{| T({\bf 0},y)-|y|| \le t \right\} \ge 1- c_1
        \cdot \exp\{-c_2 t/ \sqrt{|y|}\}.  
\]   
\end{theorem}

\subsection{A Triangle Inequality}\label{sec:triangle}

The hypothesis of Theorem \ref{Newman:Theorem1} is that the Richardson
shape $\mathcal{S}$ is \emph{uniformly curved}, that is, that there
exists $\varrho >0$ such that for each $x\in \partial \mathcal{S}$
there is a (Euclidean) ball $D^{*} (x)$ of radius $\varrho $ containing
$\mathcal{S}$ with $x$ on its surface.  Denote by $\pi
:\zz{R}^{d}\setminus \{0 \}\rightarrow \partial \mathcal{S}$ the natural projection
onto the boundary of the Richardson shape, that is, for any $x \in \mathbb{R}^{d}$
such that $x \not =0$,
\[
        \pi x =x/|x|.
\]

\begin{lemma} \label{Newman:Lemma:CircleDistance} Suppose that
$\mathcal{S}$ is uniformly curved, then there exists a constant $c>0$
such that, for all $x \not \in \mathcal{S}$ and $y \in \mathcal{S}$
\[
|x-y| \ge |x-\pi x| + c\cdot|\pi x-y|^{2}
\]
\end{lemma}

\begin{proof}
Let $|| \cdot ||_2$ be the Euclidean norm ($L_2$-norm)
on $\mathbb{R}^d$.  Any two norms on $\zz{R}^{d}$ are equivalent, and
so the Euclidean norm is equivalent to the Richardson norm: in
particular, there is a constant $\delta>1$ such that, 
for any $x \in \mathbb{R}^d$,
\[
\frac{1}{\delta}\cdot ||x||_2 \le |x| \le \delta \cdot  ||x||_2. 
\]

Let ${l}$ be the tangent hyperplane to
$D^{*}(x)$ at $\pi x$.  For any $y \in \mathcal{S}$, denote by $y_p$
the (orthogonal) projection of $y$ on ${l}$.
Since  $\mathcal{S}$ is uniformly curved (and
hence strictly convex),
\[ 
|x-y| \ge |x-\pi x| + |y-y_p|.
\]
Elementary trigonometric observations imply that, since $y \in
D^{*}(x)$ and $y_p$ is on the tangent line $l$, there exists a
constant $c_0$ that depends only on $\varrho $  and such that for all such
$y$'s
\[
||y-y_p||_2 > c_0 \cdot||\pi x-y||_2^2
\]
Hence,
\[
|y-y_p| > \frac{c_0}{ \delta^3} \cdot |\pi x-y|^2
\]
which proves the inequality  of the lemma.
\end{proof}

\section{Proof of Theorem \ref{Newman:Theorem1}} \label{sec:proof} 
\subsection{Strategy}\label{ssec:strategy} We begin by showing that it
suffices to restrict attention to a special class of initial
configurations, which we dub \emph{sliced Richardson shapes}. These
are obtained as follows: Run the Richardson model (starting from the
initial configuration with two adjacent occupied sites, one at the
origin) for a (large) time $T$, and let $Z_{T}$ be the occupied
set. Let $\xi_{T}$ be the subset of $Z_{T}$ consisting of all points with
positive first coordinates, and let $\zeta_{T}=Z_{T}\setminus \xi$.
Observe that, starting from \emph{any} nondegenerate finite initial
configuration the competition model can evolve to a sliced Richardson
shape in finite time, with positive probability. (This will occur if,
following the first time that there are adjacent Red and Blue sites,
only these sites reproduce, and only on their sides of the hyperplane
separating them.) Thus, it suffices to prove that for all sufficiently
large $T$, with positive probability the sliced Richardson shape
$\xi_{T},\zeta_{T}$ is such that
\begin{equation}\label{Newman:EquationInitialConfiguration}
        P_{\xi ,\zeta }
        \{\cap_{t\geq 0}R (t) = \emptyset \; \text{or} \;
                \cap_{t\geq 0}B (t) = \emptyset \}<4\varepsilon .
\end{equation}
Here and in the sequel $P_{\xi ,\zeta}$
will denote the probability measure governing the evolution of the
competition model under the initial condition $R (0)=\xi$, $B (0)=
\zeta$. By the Bonferroni inequality, it suffices to prove that
\begin{equation}\label{Newman:EquationInitialConfigurationRed}
        P_{\xi ,\zeta }
        \{\cap_{t\geq 0}R (t) = \emptyset  \}<2\varepsilon .
\end{equation}

The idea behind the proof of
(\ref{Newman:EquationInitialConfigurationRed}) is this: If the initial
condition $\xi ,\zeta$ is such that $\xi $ and $\zeta$ are,
approximately, the intersections of $k\mathcal{S}$ with complementary
angular sectors $\mathcal{A},\mathcal{A}^{c}$ in $\zz{R}^{d}$ based at
the origin, for large $k$, then at time $t= \delta k$ the sets $R
(t),B (t)$ should, with high probability, be approximately the
intersections of $k (1+\delta)\mathcal{S}$ with the same angular
sectors $\mathcal{A},\mathcal{A}^{c}$. This is because (1) the Shape
Theorem for the Richardson model implies that $R (t)\cup B (t)$ should
be close to $(1+\delta)k\mathcal{S}$; (2) the uniform curvature of
$\mathcal{S}$ implies that the \emph{first} occupations of vertices in
$\mathcal{A}\cap ((t+k)\mathcal{S}\setminus k\mathcal{S})$ and
$\mathcal{A}^{c}\cap ( (t+k)\mathcal{S}\setminus k \mathcal{S})$ should
(except for those near the boundaries) be by Red and Blue,
respectively; and (3) Lemma \ref{lemma:invasion} implies that, once a
region is totally occupied by Red, it must remain so (except near its
boundary) for a substantial amount of time afterward.

\subsection{Stabilization Estimate}\label{ssec:stabilization}
The key step is to show that once one of the species (say Red) has
occupied an angular sector in the  Richardson shape, it is very
unlikely for the opposite species (Blue) to make a large incursion
into this sector for some time afterward. 
Henceforth, let $d$ be the metric associated with the Richardson
norm. For any set $\zeta$ and any vertex $x$, define the distance
$d(\zeta,y)$ between $\zeta$ and $x$ to be the infimum of the distances
 $d(y,x)$ for all vertices in $y \in \zeta$.
For any point $x\in \zz{R}^{d}$ and any $r>0$, denote by
$D(x;r)$ the disk of radius $r$ centered at $x$ relative to the metric
$d$. (We shall not attempt to distinguish between open and closed
disks, as this distinction will not matter in any of the estimates.)
For $r_{1}<r_{2}$ denote by $D (x;r_{1},r_{2})$ the annular region $D
(x;r_{2})\setminus D (x;r_{1})$. For each $z\in \partial \mathcal{S}$
and any $\varrho >0$, define the 
\emph{angular sector} $\mathcal{A} (z;\varrho)$ of aperture $\varrho $
centered at $z$ by
\begin{gather*}
        \mathcal{A} (z;\varrho):= \{y \in \zz{R}^{d}\setminus \{0 \}
                \, : \, d (\pi y,z)<\varrho  \}.
\end{gather*}
Fix $r>0$, $n\geq 1$, $\beta \in (1/2,1)$ and $\alpha \in (1/2,1)$
such that $(\beta +1)/2<\alpha $, and let $A_{1}\subset
A_{2}$ be angular sectors with common center $z$ and apertures
$r<r+n^{\alpha -1}$, respectively.
 Fix $\delta \in (0,1)$, and define
\begin{align*}
        \mathcal{R}_{0}&=       \mathcal{R}_{0}^{n}=
         D (0;n/ (1+\delta),n-n^{\beta})
                \cap A_{2}\cap \zz{Z}^{d},\\
        \mathcal{B}_{0}&=       \mathcal{B}_{0}^{n}=(
                D (0;n/ (1+\delta))\cup 
                 (D (0;n/ (1+\delta), n+n^{\beta})
                \cap A_{2}^{c}))\cap \zz{Z}^{d}, \\
        \mathcal{B}_{1}&=       \mathcal{B}_{1}^{n}=(
                D (0;n)
                \cup  A_{1}^{c})\cap \zz{Z}^{d}, \\
 \mathcal{R}_{1}&=       \mathcal{R}_{1}^{n}=(
                D (0;n, n(1+\delta)-(n+\delta n)^{\beta})
                \cap  A_{1})\cap \zz{Z}^{d},
\end{align*}

\begin{lemma}\label{lemma:stabilization}
There exist constants $c_{1},c_{2}>0$ such that the following is true,
for any $n\geq 1$. If the initial configuration $\xi ,\zeta$ is such
that $\xi \supset \mathcal{R}_{0}^{n}$ and $\zeta \subset
\mathcal{B}_{0}^{n}$, then 
\begin{equation}\label{eq:stabilization}
        1-P_{\xi ,\zeta}\{
                B (\delta n)\subset \mathcal{B}_{1}^{n} \}
        \leq c_{1}n^{3d} \exp \{-c_{2} (\delta  n)^{\beta -1/2} \}.
\end{equation}
\end{lemma}

\begin{proof}
To prove (\ref{eq:stabilization})  we find exponential upper  bounds on 
\begin{gather*}\label{eq:ProbOne}
P_{\mathcal{R}_{0},\mathcal{B}_{0} } \{ B(\delta n)
 \cap \mathcal{R}_{1} \not = \emptyset \} \quad \text{and}\\
P_{\mathcal{R}_{0},\mathcal{B}_{0} } \{
% \ (B(\delta n)\cap A_1)  \cap   D(0;(n+n\delta)-(n+n\delta)^{\beta})= \emptyset 
  B(\delta n)  \cap   ( \mathcal{R}_{1} \cup \mathcal{B}_{1})^{c} 
\not = \emptyset \}.
\end{gather*}

\begin{claim}  \label{Newman:DistanceComparisonEquation}
For all  sufficiently large $n$  and for all 
$x \in \mathcal{R}_{1}$ 
\begin{equation*} 
% \label{Newman:DistanceComparisonEquation}
d(x, \mathcal{R}_{0})+ 4 (\delta n)^{\beta}   < d(x, \mathcal{B}_{0})
\end{equation*}
\end{claim}

\begin{proof}
[Proof.] 
First, observe that by Lemma \ref{Newman:Lemma:CircleDistance},
for every $x \in \mathcal{R}_{1}$ such that $|x|> {n}+n^{\beta}$, we
have
\[
d(x, \mathcal{B}_{0} ) > |x- \pi x \cdot({n}+n^{\beta})|+c \cdot n^{2\alpha-1}.
\]
Also, 
\[
d(x, \mathcal{R}_{0}) \le |x- \pi x \cdot({n}+n^{\beta})| + 2n^{\beta}. 
\]
These inequalities imply
Claim \ref{Newman:DistanceComparisonEquation} for all  $x \in 
\mathcal{R}_{1} \cap D(0; {n}+n^{\beta})^{c}$ and  all sufficiently large
$n$.

Next, if $ x \in \mathcal{R}_{1}$ is such that $n < |x| <
{n}+n^{\beta}$, then
\[ 
d(x, \mathcal{R}_{0}) \le 2n^{\beta}.
\]
Also, by Lemma \ref{Newman:Lemma:CircleDistance}, for all $y \in \mathcal{B}_{0}$,
\[
d(y, x ) > c \cdot (\delta n)^{2\alpha-1}.
\]
These two inequalities imply Claim \ref{Newman:DistanceComparisonEquation}
for all   $ x \in \mathcal{R}_{1} \cap D(0;n, n+n^{\beta})$ and all
sufficiently large $n$. 
\end{proof}

\begin{claim}\label{claim:2}
With probability $\rightarrow 1$ as $n \rightarrow \infty$, for every
$x \in \mathcal{R}_1$, every site of the ball $D(x; (\delta
n)^{\beta})$ will be colonized by time 
$\tau_x=d(x, \mathcal{R}_{0})  + 2 (\delta n)^{\beta}$. 
In particular, there exist constants $c_1>0$ and $c_2>0$ (not depending on $x$) such that 
for every $x \in \mathcal{R}_1$ 
\begin{equation*} \label{Newman:C_n:bound}
 P_{\mathcal{R}_{0},\mathcal{B}_{0} }\{ D(x, (\delta n)^{\beta}) \not \subset
 R(\tau_x)\cup B(\tau_x) \} 
\le  c_1 \cdot (\delta n)^{d \beta} \cdot \exp \{ -c_2 \cdot (\delta n)^{\beta-1/2} \}. 
\end{equation*}
\end{claim}

\begin{proof}
[Proof.]
Fix  $x \in \mathcal{R}_1$ and let
$\tau_x=d(x, \mathcal{R}_{0})  + 2 (\delta n)^{\beta}$. 
For every $z \in D(x;(\delta n)^{\beta})$,
\begin{equation} \label{Newman:C_n:distance}
d(z, \mathcal{R}_{0})+ (\delta n)^{\beta}  \le \tau_x. 
\end{equation}
Notice that 
%$d(z, \mathcal{R}_{0})  < \delta n$ and 
the number of sites in  $D(x; (\delta n)^{\beta})$ is of order  $(\delta n)^{d\beta}$.  
By  Theorem \ref{Kesten} and  (\ref{Newman:C_n:distance}), it follows 
that for some  $c_1>0$  and $c_2>0$ (not depending on $x$)
\[ 
P_{ \mathcal{R}_{0},\mathcal{B}_{0}}
\{  
D(x; (\delta  n)^{\beta})  \not \subset R(\tau_x) \cup  B(\tau_x) 
\} =
P_{\mathcal{R}_{0} \cup \mathcal{B}_{0} } \{  D(x; (\delta n)^{\beta})  \not \subset Z(\tau_x) \} \le
\]
\begin{equation*} 
\le c_1 \cdot (\delta n)^{d \beta} \exp \{ -c_2 \cdot (\delta n)^{\beta-1/2} \},
%\mbox{ for some } c_1>0 \mbox{ and }  c_2>0.
\end{equation*}
Hence, with probability $\rightarrow 1$ as $n \rightarrow \infty$, for every
$x \in \mathcal{R}_1$, every site of the ball $D(x; (\delta n)^{\beta})$ will be colonized  
at time  $\tau_x$.
 \end{proof}
%%%%%%%%%%%%%%%%%%%%%%%%%%%%%%%%%%%%%%%%%%%

 Define the boundary $\partial \mathcal{A}$ of a set $\mathcal{A} \subset \zz{Z}^{d}$
as the set of all $z \in \mathcal{A}$ that have at least one nearest neighbor that is not in $\mathcal{A}$. 
Next, for a set $\mathcal{A}   \subset \zz{Z}^d$ let $\tau(\mathcal{A})$ be  the first time at which the blue species  reaches $\mathcal{A}$.
\begin{claim}\label{claim3}
With probability $\rightarrow 1$ as $n \rightarrow \infty$, for every
$x \in \mathcal{R}_1$,  the blue species will not reach
the ball $D(x; (\delta n)^{\beta})$ by time $\tau_x$. 
In particular, there exist constants $c_1>0$ and $c_2>0$ such that  for every $x \in \mathcal{R}_1$ 
\begin{equation*} \label{Newman:D_n:bound}
P_{\mathcal{R}_{0},\mathcal{B}_{0} } \{ \tau( D(x; (\delta n)^{\beta}) ) < \tau_x \} 
\le  c_1 \cdot (\delta n)^{(d-1)(\beta+1)} \cdot \exp \{ -c_2 \cdot (\delta n)^{\beta-1/2} \}. 
\end{equation*}
\end{claim}
\begin{proof}
 Notice that
\[ 
P_{\mathcal{R}_{0},\mathcal{B}_{0} } \{ \tau( D(x; (\delta n)^{\beta})) \le \tau_x  \} \le 
  P_{\mathcal{B}_{0} } \{ Z(\tau_x) \cap D(x; (\delta n)^{\beta}) \not = \emptyset \}  
\]
\[
=P_{D(x; (\delta n)^{\beta})} \{ Z(\tau_x) \cap \mathcal{B}_{0}  \not = \emptyset \} 
\]
By Claim \ref{Newman:DistanceComparisonEquation}, for large $n$ we have  $d(x, \mathcal{B}_{0}) > \tau_x + 2(\delta n)^{\beta}$. Hence,
\[
P_{D(x; (\delta n)^{\beta})} \{ Z(\tau_x) \cap \mathcal{B}_{0}  \not = \emptyset \} \le
P_{D(x; (\delta n)^{\beta})} \{ Z(\tau_x) \not \subset D(x; \tau_x + 2(\delta n)^{ \beta })  \}. 
\]
Obviously, for every $z \in \partial D(x;(\delta n)^{\beta})$ and $y \in \partial D(x; \tau_x + 2(\delta n)^{ \beta })$,  we have 
\begin{equation*}
% \label{Newman:D_n:distance}
d(z, y )  \ge  \tau_x+ (\delta n)^{\beta}/2 .
\end{equation*}
Apply Theorem \ref{Kesten} to each pair of such vertices to get:
\[
P \{ T(z,y) < \tau_x \} < c_1 \cdot \exp \{ -c_2 \cdot (\delta n)^{\beta-1/2} \}. 
\]
The number of vertices  in $\partial D(x, (\delta n)^{\beta})$ is 
of order $ (\delta n)^{(d-1)\beta}$, and
the number of vertices  in $\partial D(x, \tau_x+2(\delta n)^{\beta})$ is
 of order at most  $(\delta n)^{(d-1)}$.
Hence,
\begin{equation*} 
%\label{Newman:D_n:Intermediate:bound}
P_{D(x; (\delta n)^{\beta})} \{ Z(\tau_x) \not \subset D(x; \tau_x + 2(\delta n)^{ \beta }) 
\le c_1 \cdot (\delta n)^{(d-1)(\beta+1)}  \exp \{ -c_2 \cdot (\delta n)^{\beta-1/2} \}. 
\end{equation*}
This finishes the proof of Claim \ref{claim3}.
\end{proof}

%%%%%%%%%%%%%%%%%%%%%%%%%%%%%%%%%%%%%%%%%%%%%%%%%%%%%%%%%%%%%%%%%%%
\begin{claim}\label{claim4}
With probability $\rightarrow 1$ as $n \rightarrow \infty$, 
there are no Blue particles in $\mathcal{R}_1$ 
at time $\delta n$.
 In particular, there exist constants $c_1>0$ and $c_2>0$ such that  
\begin{equation*}  \label{Newman:C_n+1:bound}
P_{\mathcal{R}_{0},\mathcal{B}_{0} }\{ \mathcal{R}_{1} \cap  B(\delta n) \not =
\emptyset \}  \le 
c_1 \cdot (\delta n)^{(d-1)(\beta+1)+d} \cdot \exp \{ -c_2 \cdot (\delta n)^{\beta-1/2} \}.
\end{equation*}
\end{claim}
\begin{proof}
 By Claim \ref{claim3},
for all  $x \in \mathcal{R}_{1}$ with  $\tau_x \ge \delta n$, 
\begin{equation} \label{R_1:tau_x>delta n}
P_{\mathcal{R}_{0},\mathcal{B}_{0} } \{ x  \in  B(\delta n) \} 
\le  c_1 \cdot (\delta n)^{(d-1)(\beta+1)} \cdot \exp \{ -c_2 \cdot (\delta n)^{\beta-1/2} \}. 
\end{equation}
Next, for all  $x \in \mathcal{R}_{1}$ with  $\tau_x < \delta n$, Claim \ref{claim:2} and Claim \ref{claim3}
imply that for some $c_1>0$ and $c_2>0$
%Therefore, by   (\ref{Newman:C_n:bound}) and (\ref{Newman:D_n:bound}),
\begin{equation*} \label{Newman:Ball:bound}
 P_{\mathcal{R}_{0},\mathcal{B}_{0} }\{ D(x, (\delta n)^{\beta}) \not \subset
 R(\tau_x) \} 
\le  c_1 \cdot (\delta n)^{(d-1)(\beta+1)} \cdot \exp \{ -c_2 \cdot (\delta n)^{\beta-1/2} \}. 
\end{equation*}
Hence, by Lemma   \ref{lemma:invasion}, there exist constants $c_1$ and $c_2$ such that for  every  such $x$  we have 
\begin{equation} \label{R_1:tau_x<delta n}
 P_{\mathcal{R}_{0},\mathcal{B}_{0} }\{ x  \not \in R(\delta n) \} \le 
c_1 \cdot (\delta n)^{ (d-1)(\beta+1)} \cdot \exp \{ -c_2 \cdot (\delta n)^{\beta-1/2} \}. 
\end{equation}
The number of vertices in $\mathcal{R}_{1}$ is of order $(\delta n)^d$. Thus,
by combining (\ref{R_1:tau_x>delta n}) and (\ref{R_1:tau_x<delta n}), we get
\begin{equation*} 
%\label{Newman:C_n+1:bound}
P_{\mathcal{R}_{0},\mathcal{B}_{0} }\{   B(\delta n) \cap  \mathcal{R}_{1}  \not =
\emptyset \}  \le 
c_1 \cdot (\delta n)^{(d-1)(\beta+1)+d} \cdot \exp \{ -c_2 \cdot (\delta n)^{\beta-1/2} \}.
 \end{equation*}
\end{proof}

\begin{claim} \label{claim5}
With probability $\rightarrow 1$ as $n \rightarrow \infty$, 
the blue species will not reach the set $(\mathcal{R}_{1} \cup \mathcal{B}_{1})^{c}$ 
by time $\delta n$. 
In particular, there exist constants $c_1>0$ and $c_2>0$ such that   
\begin{equation*}  \label{Newman:D_n+1:bound} 
 P_{\mathcal{R}_{0},\mathcal{B}_{0} } \{ \tau ( (\mathcal{R}_{1} \cup \mathcal{B}_{1})^{c}) \le \delta n  \}
% P_{\mathcal{R}_{0},\mathcal{B}_{0} } \{ B(\delta n)  \cap   ( \mathcal{R}_{1} \cup \mathcal{B}_{1})^{c}\not = \emptyset \}
\le  c_1 \cdot n^{2(d-1)} \cdot \exp \{ -c_2 \cdot (\delta n)^{\beta-1/2} \}. 
\end{equation*}
\end{claim}
\begin{proof}
For large $n$ the distance between the sets $\mathcal{B}_{0}$ and $(\mathcal{R}_{1} \cup \mathcal{B}_{1})^{c}$ is greater than $\delta n +(\delta n)^{\beta}$.
%Observe  that, for any  $x  \in  (\mathcal{R}_{1} \cup \mathcal{B}_{1})^{c}$,
%we have $d(x, \mathcal{B}_{0})> \delta n +(\delta n)^{\beta}$. 
The number of vertices on the boundary of $\mathcal{B}_{0} $ is  of order $n^{d-1}$.
Using the same line of argument as in the proof of Claim\ref{claim3}  we get
\[
P_{\mathcal{R}_{0},\mathcal{B}_{0} } \{  \tau ( (\mathcal{R}_{1} \cup \mathcal{B}_{1})^{c}) \le \delta n    \} \le P_{\mathcal{B}_{0} } \{  Z(\delta n)  \not \subset \mathcal{R}_{1} \cup \mathcal{B}_{1} \} \le
\]
\begin{equation*}  
c_1 \cdot n^{2(d-1)} \cdot \exp \{ -c_2 \cdot (\delta n)^{\beta-1/2} \}. 
\end{equation*}
\end{proof}
Now,   Claim \ref{claim4} and  Claim \ref{claim5}  imply 
(\ref{eq:stabilization}) and finish the proof of Lemma \ref{lemma:stabilization}:
\[
P_{\mathcal{R}_{0},\mathcal{B}_{0}} \{ B(\delta n) \not\subset \mathcal{B}_{1} \} \le
 c_1 \cdot n^{3d} \cdot \exp \{ -c_2 \cdot (\delta n) ^{\beta-1/2} \}. 
\]

\end{proof}

\subsection{Proof of
(\ref{Newman:EquationInitialConfigurationRed})}\label{ssec:proof}
Let $Z_{t}$ be the set of sites occupied by the Richardson
evolution (started from the default initial configuration) at time
$t$. Fix $T=T_{0}\geq 1$  and $\delta >0$, and set
\begin{align*}
        T_{n}&=T (1+\delta)^{n}  ,\\
        t_{n}&=T \delta (1+\delta)^{n-1}, \quad \text{and}\\
        \tau_{n}&= T_{n}-T=\sum_{j=1}^{n}t_{n}.
\end{align*}
Fix $\beta \in (1/2,1)$, and for each $n=0,1,\dotsc$ define events
\begin{align*}
       F_{n}:&=\{
         (1- T_{n}^{\beta -1})\cdot \mathcal{S}
        \subseteq 
        {\hat{Z} (T_{n})}/{T_{n}}
        \subseteq
         (1+ T_{n}^{\beta -1})\cdot \mathcal{S}\},\\
       G_{n}:&=\{
         (1- T_{n}^{\beta -1})\cdot \mathcal{S}
        \subseteq 
        {\hat{Z} (\tau_{n})}/{T_{n}}
        \subseteq
         (1+ T_{n}^{\beta -1})\cdot \mathcal{S}\}.
\end{align*}
By the Kesten-Alexander large deviation theorems (Theorem
\ref{Kesten}),
\[
         \lim_{T \rightarrow \infty }
        \sum_{n=0}^{\infty}P (F_{n}^{c})=0
\]
and so, for sufficiently large $T$, the probability is nearly $1$ that
the configuration $Z_{T}$ will be such that 
\begin{equation}\label{eq:LDSum}
        \sum_{n=0}^{\infty}P_{Z_{T}} (G_{n}^{c}) <\varepsilon .
\end{equation}

Fix an initial configuration $Z$ so that the preceding estimate holds
for $Z_{T}=Z$, and use this to construct the split Richardson shape
$\xi=\xi_{T},\zeta=\zeta_{T} $ as in section \ref{ssec:strategy}
above: $\xi $ and $\zeta$ are the subsets of $Z_{T}$ with positive and
nonpositive first coordinates, respectively. Since the union of the
Red and Blue sites in the competition model evolves as the Richardson
model, it follows from (\ref{eq:LDSum}) that if the initial
configuration is $R (0)=\xi$, $B (0)=\zeta$ then with probability in
excess of $1-\varepsilon$,
\begin{equation}\label{eq:KA}
                (1- T_{n}^{\beta -1})\cdot \mathcal{S}
        \subseteq 
        \frac{\hat{\xi } (\tau_{n})\cup \hat{\zeta} (\tau_{n})}{T_{n}}
        \subseteq
         (1+ T_{n}^{\beta -1})\cdot \mathcal{S}
\end{equation}
for all $n=0,1,2,\dotsc$. Denote by $G_{*}$ the event that
(\ref{eq:KA}) holds for all $n\geq 0$.

On the event $G_{*}$, the union of the Red and Blue regions will, at
each time $T_{n}-T$, fill a region close enough to a Richardson shape
that the estimate (\ref{eq:stabilization}) will be applicable whenever
the Red and Blue populations are restricted (at least approximately)
to angular sectors. Thus,
define sequences $A^{n}_{1}\subset A^{n}_{2}$ of
concentric angular sectors with apertures $r_{1}(n)<r_{2}(n)$
such that
\begin{gather*}
        r_{2} (n)-r_{1} (n)%=r_{3} (n)-r_{2} (n)
        =T_{n}^{\alpha -1} \quad \text{and}\\
        r_{2} (n+1)=r_{1} (n)
\end{gather*}
and with $r_{2} (0)$ chosen so that $A^{0}_{2}$ is the halfspace
consisting of all points in $\zz{R}^{d}$ with positive first
coordinates. Here $\frac{1}{2} <\beta < (\beta +1)/2<\alpha <1$ as in
Section~\ref{ssec:stabilization} above. Note that the second equality
guarantees that $A^{n+1}_{2}=A^{n}_{1}$. This in turn, together with
the fact that the sequence $T_{n}$ is increasing, implies that that
the angular sectors are nested: $A^{n+1}_{i}\subset
A^{n}_{i}$. Moreover, because $T_{n}$ is an exponentially growing
sequence and $\alpha <1$,
\[
        \lim_{n \rightarrow \infty}r_{1} (n):=r_{\infty}>0
\] 
provided $T=T_{0}$ is sufficiently large. Therefore, the intersection
\[
         A^{\infty}=\cap_{n=1}^{\infty}A^{n}_{1}
\]
is an angular sector with nonempty interior.

Finally, for each $n\geq 0$ define $H_{n}$ to be the event that at
time $\tau_{n}$ there are no blue sites in $A^{n}_{2}$ outside the
(Richardson norm) disk of radius $T_{n-1}$. (For $n=0$, set
$T_{-1}=T_{0}/ (1+\delta)$.) On the event $G_{n}\cap H_{n}$, the set
of all occupied sites is close to $T_{n}\cdot \mathcal{S}$, and the
red sites fill at least the outer layer of this set in the sector
$A^{n}_{2}$. We claim that for all sufficiently large $T$,
\begin{equation}\label{eq:final}
        P_{\xi_{T},\zeta_{T}}
        \left(\bigcap_{n=0}^{\infty} (G_{n}\cap H_{n})  \right)
        \geq 1-2\varepsilon .
\end{equation}
To see this, let $\nu $ be the smallest index $n$ such that $(G^{c}_{n}
\cup H^{c}_{n})$ occurs. Since $\nu =n$ can only occur on $G_{n-1}\cap H_{n-1}$,
\[
        P_{\xi_{T},\zeta_{T}} \{\nu =n \}
        \leq P_{\xi_{T},\zeta_{T}} (G_{n}^{c})+
         P_{\xi_{T},\zeta_{T}} (H^{c}_{n}\, | \, G_{n-1}\cap H_{n-1}).  
\]
Inequality  (\ref{eq:LDSum}) provides a bound on the sum of the first
of these terms, and Lemma \ref{lemma:stabilization} bounds the
second. Thus, for $T$ sufficiently large,
\[
        \sum_{n=0}^{\infty}P_{\xi_{T},\zeta_{T}}\{\nu =n \}< 2\varepsilon ;
\]
this proves (\ref{eq:final}).

On the event $G_{n}\cap H_{n}$, the Red species must at time
$\tau_{n}$ occupy at least the outer layer of the occupied set in the
angular sector $A^{\infty}$. Consequently, on the event $\cap_{n\geq
1} (G_{n}\cap H_{n})$, Red survives!  This proves
(\ref{Newman:EquationInitialConfigurationRed}).
\qed 

\section{Concluding Remark}\label{sec:conclusion} The preceding
argument, in addition to proving that the event that mutual survival
has positive probability, also goes part of the way towards proving
Conjecture \ref{conj}: If at a large time $T$ one of the colors (say
Red) occupies the outer layer of an angular sector, then with conditional
probability approaching $1$ as $T \rightarrow \infty$ it will occupy a
slightly smaller angular sector forever after. Since the same is true
for the other species, it follows that in at least some evolutions
Red and Blue will each occupy angular sectors. 

Unfortunately, it remains unclear what happens near the interface at
large times. Although the preceding arguments show that neither Red
nor Blue can make too deep an incursion into the other species'
sector(s), it may be possible for one to repeatedly make small
incursions across the interface that engender more (and necessarily
thinner) angular sectors in its zone of occupation. Thus, it may be
that the limit shapes exist, but consist of countably many angular
sectors. 

Finally, it remains unclear if stabilization must eventually occur on
the event of mutual survival, that is, if it is necessarily the case
that at large times $T$ the outer layer of the occupied region must
segregate into well-defined Red and Blue zones. Since local
coalescence occurs in the voter model, one naturally expects that the
same will be true in the competition model; thus far, we have been
unable to prove this.

%%%%%%%%%%%%%%%%%%%%%%%%%%%%%%%%%%%%%%%%%%%%%%%%%%%%%%%%%%%%%%%%%%%%%%

\end{document}